\documentstyle[12pt,psfig]{amsart}

\newtheorem{thm}{Theorem}
\newtheorem{lemma}[thm]{Lemma}
\newcommand{\R}{{\Bbb R}}
\newcommand{\C}{{\Bbb C}}
\newcommand{\Z}{{\Bbb Z}}
\newcommand{\Q}{{\Bbb Q}}
\newcommand{\T}{{\cal T}}
\newcommand{\X}{{\cal X}}
\newcommand{\A}{{\cal A}}
\newcommand{\I}{{\| \!|}}
\newcommand{\cO}{{\cal O}}

\begin{document}
\title{The construction of self-similar tilings}
\author{Richard Kenyon}
\address{CNRS  UMR 128, Ecole Normale Sup\'erieure de Lyon,
46, all\'ee d'Italie, 69364 Lyon, France.}
\thanks{Research at MSRI is supported in part
by NSF grant DMS-9022140.}
\date{}
\maketitle
\begin{abstract} We give a construction of a self-similar tiling of the plane
with any prescribed expansion coefficient $\lambda\in\C$ (satisfying the necessary
algebraic condition of being a complex Perron number). 

For any integer $m>1$ we show that there exists 
a self-similar tiling with $2\pi/m$-rotational symmetry group
and expansion $\lambda$ if and only if either $\lambda$ or $\lambda e^{2\pi i/m}$
is a complex Perron number for which 
$e^{2\pi i/m}$ is in $\Q[\lambda]$, respectively $Q[\lambda e^{2\pi i/m}]$. 
\end{abstract}

\section{Introduction}
By a {\bf tile} we will mean a compact subset of $\R^2$ which is the closure
of its interior.
A {\bf tiling} of the plane is a collection of tiles ${\cal T}=\{T_i\}_{i\in I}$, whose
union is $\R^2$ and which have pairwise disjoint interiors.

In \cite{Thu}, Thurston introduced the notion of {\bf self-similar tiling.}
This definition was motivated from several sources: most notably, in
tilings arising from Markov partitions for hyperbolic toral automorphisms,
one-dimensional tilings arising from substitutions, and also the Penrose
tilings, which have a simple subdivision rule to create new tilings out of old.

\begin{figure}[htb]
\centerline{\psfig{figure=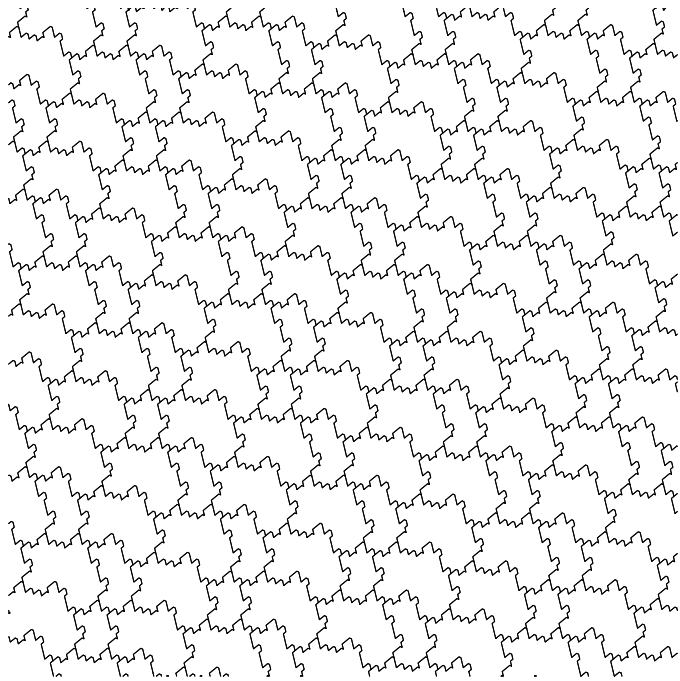,height=3.5in}}
\caption{\label{sst}Part of a self-similar tiling (see section 
\protect\ref{endos}).}
\end{figure}

A {\bf self-similar tiling} is a tiling  ${\cal T}=\{T_i\}_{i\in I}$
of the plane which has the following properties:
\begin{enumerate}
\item There is an equivalence relation $\sim$ on tiles, with a finite number of equivalence 
classes, such that $T_i\sim T_j$ implies $T_j$ is a {\it translate} of $T_i$.
\item There is a homothety $\varphi\colon\C\to\C$, $\varphi(z)=\lambda z$, such
that the image of a tile is equal to a union of tiles of ${\cal T}$,
\item If $T_i\sim T_j$, and $T_i+c=T_j$, then for each tile $T_k\subset\varphi(T_i)$,
$\varphi(c)+T_k$ is a tile of ${\cal T}$ equivalent to $T_k$.
\item The tiling is quasiperiodic (see definition below).
\end{enumerate}
An example is shown in Figure \ref{sst} (more information about this example
is given in section \ref{endos}).

Let us make a few remarks regarding this definition.

\noindent 1. The third condition simply states that two equivalent tiles ``subdivide'' in
the same way.

\noindent 2. Here by {\bf quasiperiodic} we mean, any ``arrangement'' (finite collection
of tiles, up to translation) of tiles in the tiling
reoccurs at a bounded distance from any point in the tiling.
More formally, for any $r>0$ there is an $R>0$ such that for all $x,y\in C$ 
a translate of the arrangement of tiles in $B_r(x)$ occurs as a subset of in $B_R(y)$.
This property is almost the same as the {\em local isomorphism property} 
of Radin and Wolff \cite{RW}, and is identical if one assumes in their case that
there are a finite number of local arrangements (which in our case
is a consequence of the definition).
This condition is really a compactness condition: indeed,
this condition is equivalent to the condition that the orbit closure $Y$ of the tiling
(under the action of $\R^2$ by translations) is minimal: every $\R^2$-leaf of $Y$
is dense in $Y$ \cite{RW}. 
Under reasonable assumptions on the subdivision rule for tiles,
this condition is easily checkable, see Lemma \ref{prim}, below.

\noindent 3. By replacing ``translations'' in the first and third conditions by some other
group of homeomorphisms of $\R^2$, e.g. by isometries, one arrives at a more
general definition of self-similar tiling. So technically we should refer
to our definition as that of a {\it translation}-self-similar tiling. 
However for succinctness (and to agree with Thurston's terminology) in this paper
we shall simply say
``self-similar tiling'', with the translations understood.
The results as stated below are unknown for isometry-self-similar tilings;
for similarity-self-similar tilings, see \cite{Kinfl}.
\medskip

Thurston proved that the expansion coefficient $\lambda$
of a self-similar tiling had to
be a special type of algebraic integer, a complex Perron number.
He indicated, without proof, that for each complex Perron number 
there existed a tiling. In this article we give a construction of
a self-similar tiling for any complex Perron number.

A {\bf complex Perron number} is an algebraic integer which is strictly
larger in modulus than its Galois conjugates (the other roots of its minimal
polynomial), except for $\overline\lambda$, its complex conjugate.

\begin{thm}[Thurston \cite{Thu}] The expansion of a self-similar
tiling is a complex Perron number.
\label{thth}
\end{thm}

In section \ref{const} we will prove by construction:
\begin{thm} For each complex Perron number $\lambda$
there is a self-similar tiling with expansion $\lambda$. \label{1}
\end{thm}
The methods we use are very geometric and are inspired 
by conversations with and ideas of Thurston.

This construction allows an easy corollary, whose proof we give in
section \ref{cors}:
\begin{thm}\label{rotate} $\lambda$ is the expansion of a self-similar tiling
invariant under a rotation $e^{2\pi i/n}$ ($n$ an integer) if and only if
either $\lambda\not\in\R$ and
$\lambda$ is a complex Perron number for which $e^{2\pi i/n}\in\Q[\lambda]$,
or $\lambda e^{2\pi i/n}$ is a complex Perron number and 
$e^{2\pi i/n}\in\Q[\lambda e^{2\pi i/n}]$.
\end{thm}

Note that in the case $n=2$ we conclude: any complex Perron number is the expansion
for a self-similar tiling invariant under $x\to -x$.

For an example of this theorem, consider the symmetric Penrose tiling, 
which is invariant
under rotation by $2\pi/5$. One possible expansion of this tiling
is the degree-$4$ complex Perron number
$\lambda=1+e^{2\pi i/5}$ and $\Q[\lambda]=\Q[e^{\pi i/5}]$.

In section \ref{endos} we will give a simpler construction for self-similar tilings whose
expansion coefficients are
of a particular type: those non-real complex Perron numbers $\lambda$
satisfying an equation of the form
$$\lambda^n-a\lambda^{n-1}+b\lambda+c,$$
where $a,b\geq 0,~c>1$ are integers.  This construction has the advantage of being 
very concrete.

\section{Historical remarks}

Self-similar tilings appear in several quite different contexts.
In fact up to now they have not really had a life of their own, but have rather 
been studied independently by people working in different fields.

Self-similar tilings did not really exist before the computer age, although some 
precursors can be found \cite{Golomb}. One could also trace back the roots to
the work of Adler and Weiss \cite{AW}, Bowen \cite{Bow} and Sinai \cite{Sin} 
on the construction of Markov partitions, 
which can be interpreted as the first general constructions
for self-similar tilings. 

Most of the first {\em explicit} examples of self-similar tilings were 
curiosities, found in popular books on fractals with
fanciful names such as the ``Heighway dragon curve'' and 
``Gosper flowsnake''. 

The Penrose tilings \cite{Pen} were another curiosity until some 
nonperiodic crystals displaying $5$-fold ``symmetry''
were actually discovered in nature; these ``quasicrystals"
gave life to a new field of crystallography.

The theory of wavelets is another new field which has
caused people to become interested
in ``self-replicating'' tilings (self-similar tilings with one tile type), 
and much work has been done since on this special class of self-similar tilings \cite{Ban,
GH,St,LW,Vince}.

Our current construction has its origins in another point of view.
Working in the field of symbolic dynamics, 
Lind \cite{Lind} gave a characterization of the largest eigenvalues of nonnegative integer
matrices, which gives as corollary a 
characterization of the expansions of self-similar
tilings of the line. 

Then Thurston \cite{Thu} in 1989 took a step up from this result,
defining self-similar tilings and proving the neat Theorem \ref{thth} above. 
As a thesis student of Thurston, I extended his result to 
$n$-dimensional self-affine tilings (in which the expansion is a linear map of
$\R^n$ diagonalizable over $\C$) \cite{Ken}.

Not much later Praggastis \cite{Prag} gave a more explicit construction
of Markov partitions for (certain) hyperbolic toral automorphisms using the
theory of self-similar tilings. 
This work is also very close to that of Vershik \cite{Ver}
who constructs such Markov partitions {\it arithmetically}, that is, 
using radix representations.

Self-similar subdivisions (of various sorts) occur nowadays in the theory of rational maps,
Kleinian and hyperbolic groups (via the action on the sphere at infinity or space of
geodesics using an automatic structure), wavelets, crystallography (in 
quasicrystals), symbolic dynamics and number theory (in radix representations), and
logic (in Wang tiles, aperiodic tilings, and computability: surprisingly,
the only known aperiodic tilings are built from self-similar patterns).

\section{Background}

\subsection{the subdivision rule}
The {\bf tile types} in a self-similar tiling are the distinct
equivalence classes of tiles.
To a self-similar tiling with $n$ tile types is associated a nonnegative $n\times n$
integer matrix $M=(m_{ij})$, the {\bf subdivision matrix},
in which $m_{ij}$ is the number of tiles of type $j$
that a tile of type $i$ subdivides into upon scaling by $\lambda$. By property $3$
this matrix is well-defined.

The Perron (i.e. largest) eigenvalue 
of this matrix is $\lambda\overline\lambda$, whose corresponding 
eigenvector is the vector of tile areas. 

Recall that a nonnegative matrix is called {\bf primitive}
if some power is strictly positive. In the case of a self-similar tiling,
the subdivision matrix must be primitive, since every tile, when scaled by $\lambda^n$
for large enough $n$, must contain a copy of every tile type by quasiperiodicity.
Primitivity also implies that the spectral radius 
$\lambda\overline{\lambda}$ is a {\bf Perron number},
that is, $\lambda\overline{\lambda}$
is a real algebraic integer strictly larger than the modulus of any Galois conjugate.

We have the following result (whose proof is an application of the 
definitions) which we shall use later.

\begin{lemma}[\cite{Prag}, see also \cite{Kinfl}] A tiling which satisfies the first three 
hypotheses of a self-similar tiling, has a tile
with the origin in its interior and has primitive subdivision matrix, is
quasiperiodic, and hence self-similar.
\label{prim}
\end{lemma}

\subsection{The Delauney triangulation}
We review here the definition of the Delauney triangulation.
Given a discrete set of points $P$ in $\R^2$, the Delauney triangulation
is the triangulation of the convex hull of $P$, with triangles having
vertices in $P$, with the property that, for any triangle $t=(v_1,v_2,v_3)$,
the circumcircle of $t$ contains no point of $P$. Such a triangulation
exists (for example, one can see this by using stereographic projection 
of the plane to the sphere in
$\R^3$; the convex hull in $\R^3$ of the image of $P$ is a polyhedron whose
faces are (generically) the triangles in the Delauney triangulation) 
and is unique unless some four or more points lie on the same
empty circle. 

Another property that we will use is the local nature of the triangulation:
if for some small $\epsilon$ 
there is a point of $P$ within a distance $\epsilon$ of any point of $B_1(x)$, 
then the set of triangles with vertex $x$ depend only on those points in $B_{3\epsilon}(x)$. 
(To see this, notice that an empty circle passing through $x$ has
diameter at most $2\epsilon$.)

\section{Proof of Theorem \ref{1}}
\label{const}
\subsection{sketch}
The idea of the construction is as follows. We define
the set $\T$ of ``archtypes'' of the tiles in our eventual tiling
to be the set of all triangles
with vertices in an appropriate lattice and edges of bounded length.
When we multiply a triangle in $\T$ by $\lambda^n$ for large enough $n$,
we can subdivide it approximately into triangles in $\T$.
(This is where we use the complex-Perron-number property of $\lambda$). 
However we need a way to 
define the subdivision of triangles so that,
if two triangles are adjacent, their subdivisions agree at their common boundary.
So rather than just keep track of a triangle we keep track of a
triangle and all its immediate neighbors. This gives a new larger
set of archtypes $\X$. The subdivisions of elements of $\X$ can 
now be defined so that they {\it agree} whenever there is an overlap between
them.

One small hitch is to make this subdivision rule primitive.
To ensure this we define a special tile $T_0$ that occurs in the subdivision
of any tile, including itself. Our set of tiles will then be
the strongly connected
component of the ``subdivision graph'' which contains $T_0$.

Finally we redraw the boundaries of the elements of $\X\cup\{T_0\}$ by a
recursive process so that
the subdivision is exact, not just approximate.

\subsection{The set-up}
Let $\lambda$ be a complex Perron number of degree $d$. Let $q[x]\in\Z[x]$ be
the minimal
polynomial for $\lambda$; it is monic and of degree $d$.
Furthermore the roots $\lambda_0,\lambda_1,\ldots,\lambda_{d-1}$
of $q(x)$ satisfy $|\lambda_i|<|\lambda|$
except for the roots $\lambda$ and $\overline{\lambda}$.
The roots of $q(x)$ are called the {\bf Galois conjugates} of $\lambda$.
Assume $\lambda_0=\lambda$.

Let $K=\Q[\lambda]$. 
For each $i$ with $0\leq i\leq d-1$ there is an
embedding $\sigma_i\colon K\to\C$, 
sending an element $p(\lambda)$ (a polynomial in $\lambda$ with rational
coefficients) to $p(\lambda_i)$. If $\lambda_i$ is real the image of 
this embedding is contained in $\R$. 

Let $r$ be the number of real Galois conjugates, and $2c$ the number
of non-real Galois conjugates of $\lambda$, so that $r+2c=d$.

There is a natural embedding $\sigma\colon K\to\R^r\times\C^c\cong\R^d$ 
which is the product
of these individual embeddings $\sigma_i$ (taking only one embedding for each
complex conjugate pair). 

At this point we must distinguish between the case $\lambda$ real and the case $\lambda$ 
is nonreal.

If $\lambda$ is nonreal, let $W=\R^d$; let $\pi$ be the linear projection from $W$
to $\C$ such that $\forall x\in K$, $\pi\sigma(x)=x$. Let $m_\lambda$ be the linear map
$\sigma(K)\to\sigma(K)$ induced by multiplication by $\lambda$ in $K$, i.e.
$m_\lambda(\sigma(x))=\sigma(\lambda x)$.
Then $m_\lambda$ extends to a linear map of $\R^d$ whose eigenvalues
are exactly the $\lambda_i$, and eigenspaces are given by the coordinate
axes/planes.
Let $\A=\cO$ denote the ring of algebraic integers in $K$. 
The image $\sigma(\A)$ is a discrete lattice in $W$ invariant under $m_\lambda$.

If $\lambda$ is real let $W=\R^d\times\R^d$ and $\pi$ be the projection
$\pi\colon W\to\R^2$ satisfying $\pi(\sigma(x)\times\sigma(y))=(x,y)$.
Let $m_\lambda$ be the linear map on $W$ induced by multiplication by
$\lambda$ on $K$. The eigenvalues of $m_\lambda$ are the $\lambda_i$,
each with multiplicity $2$.
Let $\A=\cO\times\cO$ where $\cO$ is the ring of algebraic integers in $K$. 
The image $\sigma(\A)$ is again discrete lattice in $W$ invariant under $m_\lambda$.

In either case let $V_\lambda\subset W$ be the eigenplane (2-dimensional eigenspace)
for $m_\lambda$
corresponding to the eigenvalue $\lambda$. Then the map $\pi\colon W\to\R^2$ 
can be thought of as projection along eigenvectors of $W$ onto $V_\lambda$.

We use the Euclidean metric on both $\R^2$ and $W$. 
If $x\in \R^2$ or $x\in W$ let $N_R(x)$ 
denote the neighborhood of radius $R$ around $x$ in $\R^2$ or
$W$, respectively.
Similarly for a set $S\subset \R^2$ or $S\subset W$ 
let $N_R(S)$ denote the $R$-neighborhood of $S$.
For $v\in K$ let $\I v\I$ denote the distance from $\sigma(v)$ 
to the plane $V_\lambda$. This is the length of the ``vertical'' component of 
$\sigma(v)$.

For a triangle $t\subset\R^2$ with vertices $v_1,v_2,v_3\in K$,
let $\sigma(t)$ denote the triangle in $W$ with vertices $\sigma(v_1),
\sigma(v_2),\sigma(v_3)$.

\begin{lemma} Given a triangle $t\subset\R^2$ with vertices $v_1,v_2,v_3\in\A$
and an $\epsilon>0$,
there is an integer $n$ such that $\sigma(\lambda^nt)$ is almost 
parallel to $V_\lambda$, that is, the projection $\pi$ 
is a $(1+\epsilon)$-biLipschitz mapping from $\sigma(\lambda^nt)$ to
its image $\lambda^nt$.
\label{eflat}
\end{lemma}

\noindent{\bf Proof.}
A non-zero vector $v\in W$ lying in the 
plane of $\sigma(t)$ has a non-zero component in the direction of the plane $V_\lambda$.
Since $\lambda$ is the strictly largest eigenvalue in modulus of $m_\lambda$,
for $n$ sufficiently large the component in direction $V_\lambda$ of
$m_\lambda^nv$ dominates all the others, and so $m_\lambda^n(v)$ 
lies almost parallel to $V_\lambda$.  The result follows.
\hfill{$\Box$}

\subsection{The actors.}
\begin{figure}[htbp]
\centerline{\psfig{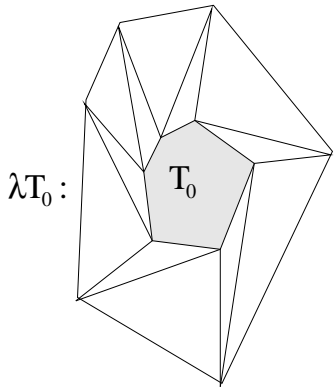}}
\caption{\label{T0}The triangulation of $\lambda T_0$.}
\end{figure}

Let $T_0\subset\R^2$ be a polygon with vertices in $\A$, having the origin
in its interior and sufficiently ``round'' so that $\lambda T_0$ contains
$T_0$ in its interior. Triangulate the annulus $\lambda T_0-T_0$ with triangles
whose vertices are in $\A$ (see Figure \ref{T0}).

Let $M$ be an integer large enough so that:
\begin{enumerate}
\item $M>2|\lambda|^2\I v_1-v_2\I $ for any two vertices $v_1,v_2$ in the triangulation
of $\lambda T_0-T_0$.
\item $M>|\lambda|^2|v_1-v_2|$ for any two such vertices.
\item For any $x\in W$, $N_{M/(2|\lambda|)}(x)$ contains a point of  the discrete lattice
$\sigma(\A)$.
\end{enumerate}

Let $\T$ be the set 
of triangles $t=(v_1,v_2,v_3)$ with noncollinear vertices $v_1,v_2,v_3\in \A$
satisfying $\I v_i-v_j\I <3M$ and $|v_i-v_j|<3M$ for all $i,j\in\{1,2,3\}$.
For any $t\in\T$, $\sigma(t)$ is a triangle with vertices in the discrete lattice 
$\sigma(\A)$ and edges of bounded length. 
Thus $\T$ is a finite set up to translation, i.e. if we consider 
$(v_1,v_2,v_3)=(v_1+a,v_2+a,v_3+a)$ for all $v_1,v_2,v_3,a\in\A$.

\begin{figure}[htb]
\centerline{\psfig{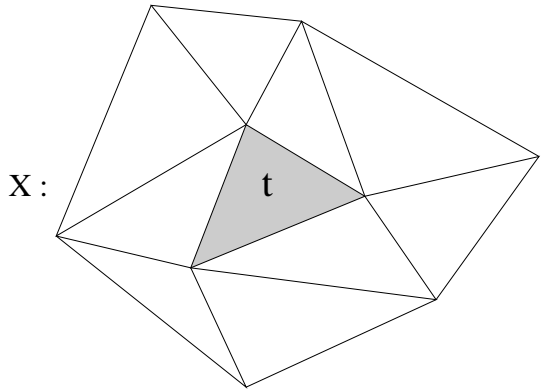}}
\caption{\label{surr}A surrounding.}
\end{figure}

Let $\X$ be the set of {\bf surroundings} of elements of $\T$ by elements of
$\T$. By definition a surrounding of $t\in \T$ is a collection $X$ of elements
of $\T$, with $t\in X$,
which tiles in an edge-to-edge fashion a neighborhood of $t$ and such 
that every element of $X$ touches $t$ in at least a vertex (see Figure 
\ref{surr}.) Note that $\X$ is also a finite set up to translation.

Recall that for $t=(v_1,v_2,v_3)\in\T$ by $\sigma(t)$ we mean 
the triangle in $W$ with vertices $(\sigma(v_1),\sigma(v_2),\sigma(v_3))$.
For $X\in\X$, by $\sigma(X)$ we mean the union of the $\sigma(t)$'s 
for each triangle $t$ in $X$. So $\sigma(X)$ is a piecewise flat surface in $W$.

\begin{figure}[htb]
\centerline{\psfig{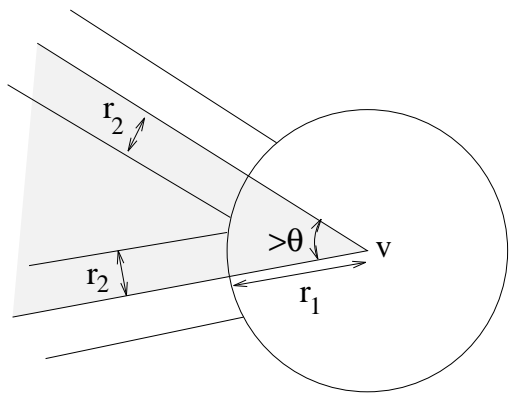}}
\caption{\label{nbhds}Neighborhood of a vertex.}
\end{figure}

Let $\theta>0$ be the smallest vertex angle of any triangle in $\T$.
Let $r_2=2M$ and $r_1=r_2/sin(\theta/2)$. Then in any homothetically scaled
copy of a triangle in $\T$, the $r_2$-neighborhoods of any two edges
do not intersect outside of the $r_1$-neighborhoods of the vertices (Figure
\ref{nbhds}).

Let $n$ be an integer large enough so that for each $t\in\T$:
\begin{enumerate}\item
$\lambda^nt$ is $1/2$-flat in the sense of Lemma
\ref{eflat} (the projection $\pi$ from $\sigma(\lambda^nt)$ is $3/2$-biLipschitz).
\item
the inscribed radius of $\lambda^nt$ is at least $2{\rm diam}(\lambda T_0)+2r_2$.
\end{enumerate}

\subsection{subdividing an element of $\X$}
\label{subd}
Let $X\in\X$ be a surrounding of $t\in\T$. 
Each triangle in $\lambda^nX$ is large and flat by our
choice of $n$. The subdivision proceeds in three stages:
first we must subdivide a neighborhood of the vertices of $\lambda^nt$, then a neighborhood
of the edges of $\lambda^nt$, and then $\lambda^nt$ itself.
(In the final construction we'll use $n+2$ in place of $n$; the rules are
the same, though.)

\subsubsection{near a vertex}
\label{subverts}
Let $X\in\X$ be a surrounding of $t\in\T$.
Let $v$ be a vertex of $t$. Let $N=N_{r_1}(\lambda^nv)\subset\lambda^nX$.
Let 
$$Y_v=\{y\in\A\cap N: \I y-\lambda^nX\I <M/|\lambda|\}.$$ 
Here by $\I y-\lambda^nX\I $ we mean
the distance vertically from $\sigma(y)$ to $\sigma(\lambda^nX)$. This is the Euclidean 
distance $d(\sigma(y),z)$ where $z\in W$ is defined by $\pi(z)=y$ 
and $z\in\sigma(\lambda^nX)$.  
Then $Y_v$ is a discrete subset of 
$N$. (Later on $Y_v$ will be the vertices in a triangulation of $N$ in
which the triangles are in $\T$.)

\subsubsection{near an edge}
\label{subedges}
Let $e=e_{ij}=\overline{v_iv_j}$ be an edge of $t$ and let
$$N'=N_{r_2}(\lambda^ne)-(N_{r_1}(\lambda^nv_i)\cup N_{r_1}(\lambda^nv_j))$$
be the $r_2$-neighborhood of $\lambda^ne$
which is outside the vertex-neighborhoods defined
in the previous paragraph. Let 
$$Y_e=\{y\in\A\cap N': \I y-\lambda^nX\I <M/|\lambda|\}, $$
where $\I y-\lambda^nX\I$ is as defined above.
Again, $Y_e$ is a discrete subset of $N'$. 

\subsubsection{for a triangle.}
Once the sets $Y_v,Y_e$ have been defined as above for each vertex and edge
of $t$, define 
$$Y_t=\{y\in\A\cap \lambda^nt: \I y-\lambda^nX\I <M/|\lambda|\}.$$
Let $$Y=Y_t\cup Y_{v_1}\cup Y_{v_2}\cup Y_{v_3}\cup Y_{e_{12}}\cup Y_{e_{23}}\cup
Y_{e_{31}};$$ then $Y$ is a discrete set of points in the $r_2$-neighborhood
of $\lambda^nt$.

Triangulate the set $Y$ using the Delauney triangulation, leaving out
triangles whose circumcircle is not completely contained in $N_{r_2}(\lambda^nt)$.
Each triangle $t'=(w_1,w_2,w_3)$ in this triangulation satisfies: 
\begin{equation}\label{abssmall}
|w_i-w_j|<M/|\lambda|\end{equation}
since the circumcircle of $t'$ contains no point of $Y$, and so must have diameter
bounded by $M/|\lambda|$ by condition (3) on $M$. Furthermore
each triangle satisfies  
\begin{equation}\label{htsmall}
\I w_i-w_j\I <\I w_i-\lambda^nX\I +\I \lambda^nX-w_j\I +
\frac12|w_i-w_j|< \frac M{|\lambda|}(1+1+\frac12)
\end{equation}
here using the fact that $\sigma(\lambda^nX)$ is $\frac12$-flat (condition (1) on $n$).

Since $|\lambda|>1$, these imply $t'\in\T$. We have defined a subdivision of 
a neighborhood of $t$ into triangles in $\T$.
In addition, for each such triangle $t'$, we have $\lambda t'\in\T$
by (\ref{abssmall}) and (\ref{htsmall}).

If two surroundings $X_1,X_2$ (of $t_1,t_2$ respectively)
overlap, i.e. $t_2\in X_1$ and $t_1\in X_2$, then
their subdivisions agree near where $\lambda^nt_1$ intersects
$\lambda^nt_2$ (by the local nature of the Delauney triangulation), 
so that $X_1\cup X_2$ has a well-defined subdivision.

\subsubsection{the central tile}
\label{T0inthere}
We need to alter slightly this subdivision rule for triangles so as
to insert a copy of the region $\lambda T_0$ in the subdivision. 
Let $t$ and $X$ be as before with $Y_v,Y_e$ defined for each vertex and edge.
Put a translate $C$ of $\lambda T_0$ in the interior of 
$\lambda^nt$ so that for any vertex $w$ 
of $C$, we have 
\begin{equation}\label{condC}
\I w-\lambda^nt\I <M,
\end{equation}
and $C$ does not intersect
the sets $N_{r_2}(\lambda^ne)$ for edges $e$ of $t$. This can be
accomplished by conditions (1) on $M$ and (2) on $n$: just put $C$ near the
center of the inscribed circle of $\lambda^nt$, so that vertices
of $C$ are in $\A$, and so that $\sigma(C)$ is (vertically)
within $M$ of $\sigma(\lambda^nt)$.
Take $$Y'=\{y\in Y\mid y\not\in{\rm interior}(C)\}.$$

As before use the Delauney triangulation of $Y'$, except near $C$: 
the edges of the boundary of $C$ may unfortunately 
not be edges of Delauney triangles.
Take only those Delauney triangles whose circumcircle doesn't penetrate into $C$.
The remaining untriangulated area is an annulus around $C$, and for any vertex $x$
of the outer boundary $Z$ of this annulus, the distance to the inner boundary $\partial C$ 
is at most $M/|\lambda|$,
since $x$ is on a circle of diameter at most $M/|\lambda|$ which touches $C$.
In particular $x$ is at distance at most $2M/|\lambda|$ from a vertex of $\partial C$,
since the edges of $\partial C$ have length at most $M/|\lambda|^2<M/|\lambda|$.
So we can triangulate this annulus in the following manner: for each edge $c_1c_2$ of
$\partial C$, connect $c_1,c_2$ to their closest common neighbor on $Z$.
Connect the remaining vertices $c_3$ on $Z$ to the only place now possible on $\partial C$
which doesn't cross edges already present.
Each such triangle has edge lengths at most $3M/|\lambda|$.
By the condition (\ref{condC}) on $\T$, each triangle is in $\T$.

\begin{figure}[htbp]
\centerline{\psfig{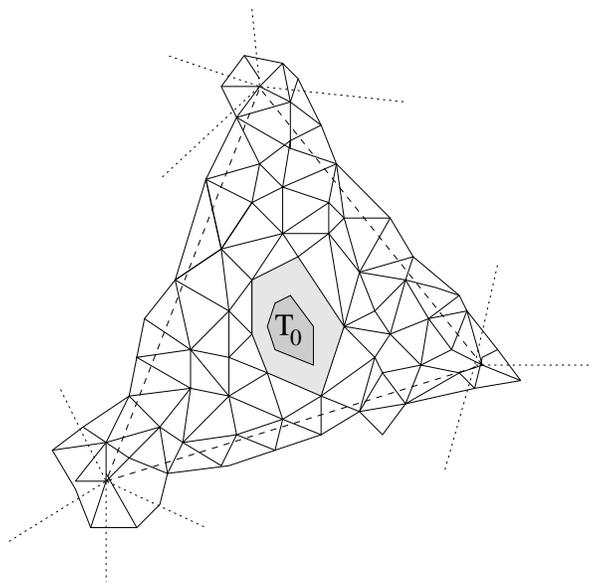}}
\caption{\label{trisubd}Subdividing a triangle.}
\end{figure}

This defines a subdivision
of a neighborhood of $t$ into triangles in $\T$ and a region $C$
which is tiled by a translate of the original tiling of $\lambda T_0$.
Furthermore triangles $t$ near the edges of $\lambda^nt$ satisfy: $\lambda t\in \T$.
(Figure \ref{trisubd}).

\subsection{the tiling}
We can now define the tiling of the plane as follows.
Start with a tile $T_0$ at the origin.
When we multiply by $\lambda$, the image $\lambda T_0$ subdivides
into $T_0$ union the given triangulation of $\lambda T_0-T_0$.
We multiply by $\lambda$ $n+2$ more times, mapping the annuli $\lambda T_0-T_0$ 
successively to $\lambda^2T_0-\lambda T_0$, $\lambda^3T_0-\lambda^2T_0$,
and so on, up to $\lambda^{n+3}T_0-\lambda^{n+2}T_0$. At this stage we
have a concentric sequence of annuli (each with a triangulation homothetic to the
triangulation in $\lambda T_0$), and a tile $T_0$
at the origin.

Each triangle except those in the innermost and outermost annulus has
a surrounding by triangles. The surroundings of triangles in the second innermost
annulus are in $\X$, since the triangles in the first three annuli are in $\T$
by conditions (1) and (2) on $M$. 
So each surrounding is a homothetic copy of an element of $\X$.
At this point the surroundings of
the triangles in the 2nd outermost
annulus are of the form $\lambda^n X$ for some $X\in \X$. 
Now subdivide (as in section \ref{subverts},\ref{subedges})
all vertices and edges which are on the boundary between
the outermost and 2nd outermost annulus (i.e. on $\lambda^{n+2}\partial T_0$). 

Multiply once more by $\lambda$. Subdivide all triangles in the 2nd outermost
annulus $\lambda^{n+3}T_0-\lambda^{n+2}T_0$ 
using the subdivision as in sections \ref{subverts},\ref{subedges}, \ref{T0inthere}. 
That is, first subdivide any vertices which have not yet been subdivided (i.e. were
not subdivided at the previous stage), then any edges which have not yet been subdivided,
then finally subdivide triangles when all their vertices and edges have been subdivided.
Those triangles adjacent to $\lambda^{n+2}\partial T_0$ 
have a previously defined subdivision
on one or more of their vertices and edges. The new subdivisions 
will by construction agree with the old
subdivisions near these vertices and edges.
Note that in this way the subdivision of a triangle depends only on the 
triangles in its surrounding (i.e. on the appropriate element of $X$).

We now continue multiplying by $\lambda$ and subdividing, with the following rule:
label the triangles by their ``size'':
those in $\lambda T_0-T_0$ are labeled $1$, the image
of a triangle of label $k$ is a triangle of label $k+1$ unless $k=n+2$, 
in which case the image is subdivided into a copy of $T_0$ and lots of small
triangles in $\T$ of label $1$ again. 

As we continue multiplying these
tilings fill out the whole plane, and there is a unique limiting tiling 
since larger and larger portions around the origin are fixed.

Adjacent triangles in the limiting tiling
have the same size or else their size differs by $1\bmod n+2$. 
If adjacent triangles have size both $<n+2$, then they don't subdivide
on multiplication by $\lambda$.
If the sizes are $n+2,n+2$ then their subdivisions agree by definition.
If the sizes are $n+1,n+2$, then the $n+2$-sized tile subdivides, so the
common egde or vertex subdivides. At the next multiplication by $\lambda$,
the other tile will subdivide with a subdivision which by construction matches the 
existing subdivision along the common edge or vertex. 
If adjacent sizes are $1$ and $n+2$, then they just came from a pair of tiles
of sizes $n+2$ and $n+1$, so again their subdivisions agree.

\subsection{refining the boundaries}

We have defined a subdivision rule for triangles, or rather, elements of $\X$,
and a tiling of the plane with scaled elements of $\X\cup\{\lambda T_0\}$.
To make this a true self-similar tiling we need to redraw the boundaries
of the triangles and of $T_0$ so that the image of a tile covers {\em exactly}
a set of tiles in the tiling.

Let $X\in\X$ be a surrounding of $t\in \T$. 
For each edge $e$ in $t$ (between vertices $v_i$ and $v_j$ of $t$)
draw a polygonal arc $\alpha_e$ from 
$\lambda^{n+2}v_i$ to $\lambda^{n+2}v_j$, which is contained in 
$N_{r_2}(\lambda^{n+2}e)$, and has edges contained in the set of edges
of triangles in the triangulation of $N_{r_2}(\lambda^{n+2}t)$. 

The part of the path $\alpha_e$
in $N_{r_1}(\lambda^{n+2}v_i)$ should depend only on
this neighborhood; similarly the part of the path in $N_{r_1}(\lambda^{n+2}v_j)$
should depend only on that neighborhood, and the part of the path along
the edge should depend only on $N_{r_2}(\lambda^{n+2}e)$. 
This is to guarantee that the paths for different surroundings agree whenever
there is overlap. For example, if triangles $t_1$ and $t_2$ are adjacent
along an edge $e$, then the path $\alpha_e$
along $\lambda^{n+2}e$ will be the same 
for the surrounding of $t_1$ as for $t_2$.

\begin{figure}[htb]
\centerline{\psfig{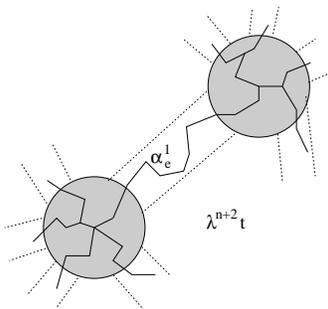}}
\caption{\label{pathchoice}Choosing a path around $t$.}
\end{figure}

Furthermore the choice of paths in vertex neighborhoods $N_{r_1}(\lambda^{n+2}v_i)$ should be
made so that the paths going out along different edges incident at the vertex
should be {\bf non-crossing} and {\bf diverging} (although such paths may have common initial
segments, they are disjoint except for a connected initial segment).
See Figure \ref{pathchoice}.

Define $\beta_e=\lambda^{-n-2}\alpha_e$, so that $\beta_e$
is a polygonal arc between vertices of $t$ which follows closely the edge $e$.
Define $\gamma^1_t$ to be the concatenation of the 3 edge paths 
$\beta_{e_{12}},\beta_{e_{23}},\beta_{e_{31}}$, removing any
``backtracking'' 
which may occur near the three vertices: that is, if the last
edge of $\beta_{e_{12}}$ is the same as the first edge of $\beta_{e_{23}}$ in the
reverse direction, remove these two edges from $\gamma^1_t$, and so on. 
So $\gamma^1_t$ is a simple closed polygonal path
running once around $t$. 

The paths $\lambda^{n+2}\gamma_t^1$ 
are composed of edges of triangles in $\T$. 
To each such edge $e_j$ is associated its arc $\beta_{e_j}$.
Define $\gamma_t^2$ by replacing each edge $\lambda^{-n-2}e_j$ in $\gamma_t^1$
with a copy of $\lambda^{-n-2}\beta_{e_j}$, again removing backtracks near the vertices
of the $e_j$. Inductively define
$\gamma_t^{k+1}$ by replacing each edge $\lambda^{-(n+2)k}e_j$
in $\gamma_t^k$ by a copy of $\lambda^{-(n+2)k}\beta_{e_j}$ and removing backtracks.

These paths $\gamma_t^k$ converge as $k\to\infty$ to a simple closed curve $\gamma_t$
running around $t$:
$\gamma_t^k$ lies in the $c\lambda^{-k(n+2)}$-neighborhood 
of $\gamma_t^{k-1}$, for some
constant $c$,
so the paths converge in the Hausdorff metric. The limits are arcs by our choice
of $n$, which is so large that any potential non-injectivity in $\gamma_t^k$ must
occur near a vertex of $\gamma_t^{k-1}$, where by construction it cannot occur.

For each $X\in\X$, the Jordan curve $\gamma_t$ forms the boundary of a tile $T_X$. 
A similar construction gives a Jordan curve for the boundary of $T_0$.
These tiles $T_X$ have the property that $\lambda^{n+2}T_X$ 
subdivides exactly into other tiles in $\{T_0\}\cup \{T_{X'}:X\in\X\}$.
Also $\lambda T_0$ subdivides exactly into a copy of $T_0$ and certain $T_X$.

\subsection{proof of self-similarity}

Let $S=\{T_0\}\cup S_1\cup S_2\cup \ldots\cup S_{n+2}$, where 
$S_i=\{\lambda^iT_X\mid X\in\X\}$.
We have now constructed a tiling of the plane, with tiles which are 
translates of tiles in a subset of $S$ (not every element of $S$ may actually occur
in the tiling). By construction when we multiply 
the tiling by $\lambda$ and subdivide each tile according to 
the rules established above ($T_0$ subdivides into the tiles as defined
in section \ref{T0}, each tile in $S_i$ maps homothetically to a tile
in $S_{i+1}$ unless $i=n+2$ in which case a tile subdivides as
described in section \ref{T0inthere}).

Furthermore tiles of the same type (i.e. those arising from the same
element of $X$, and which are of the same size) subdivide in the same way.

Furthermore each tile occurring in the tiling arises from eventual subdivision of
the tile $T_0$, since $T_0$ has the origin in its interior.
Each tile eventually has a copy of $T_0$ in its subdivision. 
So the subdivision matrix is primitive, and hence the tiling
is quasiperiodic (Lemma \ref{prim}).

This completes the proof.

\section{Tilings invariant under rotations}
\label{cors}
Let $T$ be a self-similar tiling 
invariant under $e^{i\theta}$ and having expansion constant $\lambda$.
If $\lambda$ is real, then by redefining
the subdivision rule we have that $T$ is self-similar with expansion 
$e^{i\theta}\lambda$. Thus we can assume $\lambda$ is not real unless $e^{i\theta}=-1$.

In \cite{Kinfl} it is shown that for any
self-similar tiling with expansion $\lambda\not\in\R$ there is a homothetic
copy $K'$ of $\Q[\lambda]$ such that translations between tiles
of the same type must be in $K'$. So if the tiling is invariant
under $x\to xe^{i\theta}$, then if $a\not=0$ is a translation between tiles
of the same type, both $a$ and $ae^{i\theta}$ are in $K'$, so
their ratio $e^{i\theta}$ must be in $\Q[\lambda]$.

In case $e^{i\theta}=-1$ we also have $e^{i\theta}\in\Q[\lambda]$.
\medskip

To construct a tiling invariant under rotation by $\theta=2\pi/m$, modify
the construction of section \ref{const} as follows.

Select the tile $T_0$ to be invariant under $e^{i\theta}$. Then the
region between $\lambda T_0$ and $T_0$ has a $e^{i\theta}$-invariant triangulation.

The construction now proceeds as before. Note that $e^{i\theta}$ is an algebraic
integer, and multiplication by $e^{i\theta}$
acts as an {\em isometry} on the space $W$ 
and the lattice $\sigma(\A)$.
Since the subdivision defined in section \ref{subd} 
depends solely on the metric properties of the 
objects involved, the subdivision is natural in the sense that two elements 
$X_1,X_2\in \X$ (surroundings of $t_1,t_2$ respectively)
which are isometric copies of each other have subdivisions which are
isometric. (Note that even if $X_1$ and $X_2$ overlap,
the naturality of the subdivision implies that they will subdivide
isometrically). 

The only place in the construction where we again had a choice
was in defining the paths running around the tiles. However we claim
that we can make this choice in a way which only depends on the isometry
type of the relevant vertex or edge. The only difficulty arises when
isometric surroundings overlap. If
$t_1$ and $t_2$ meet at a vertex $v$ and the isometry from $t_1$ to $t_2$ takes 
$v$ to $v$, then the subdivision near $v$ is already 
$\theta$-invariant, so one can choose paths running out from that vertex
in a $\theta$-invariant fashion. Similarly if $t_1$ and $t_2$ meet along an edge
$e$, and the isometry taking $t_1$ to $t_2$ takes $e$ to $e$,
then the subdivision of $e$ is invariant under this isometry
and so one can choose the path $\alpha_e$ to also be invariant.
These choices imply the naturality of the entire construction.

So the resulting tiling will be invariant under $e^{i\theta}$.

\section{Construction from free group endomorphisms}
\label{endos}
In this section we give a construction for self-similar tilings whose
expansion coefficient is a complex Perron number which satisfies 
\begin{equation}\label{pqr}\lambda^n-p\lambda^{n-1}+q\lambda+r=0
\end{equation}
for some $n\geq 3$ and integers $p,q\geq 0,~r\geq 1$.

We will first do the case $n=3$.
Let $a,b,c$ be vectors pointing in different directions in $\R^2$. Let
$F$ be the set of compact polygonal paths starting at the origin,
each of whose edges is a translate of $\pm a,\pm b,$ or $\pm c$, 
and which are non-backtracking in the sense that along the path
a segment $x$ is never immediately followed by $-x$. 

Thus to each element in $F$ is associated a unique element in the
free group on three symbols $F(a,b,c)$. Indeed, 
one can easily define a product on $F$ in such a way
that the map $F(a,b,c)\to F$ is a group isomorphism: to obtain the product
of two paths, translate the second to the end of the first, and then
cancel any ``backtrackings''.

Let $f\colon F(a,b,c)\to F$ be this isomorphism.

Define an endomorphism $\phi\colon F(a,b,c)\to F(a,b,c)$ by:
\begin{eqnarray*}\phi(a)&=&b\\
\phi(b)&=&c,\\
\phi(c)&=&c^pa^{-r}b^{-q}.
\end{eqnarray*}

Consider the three commutators $[a,b]=aba^{-1}b^{-1},[b,c],$ and $[a,c]$; 
they represent three closed paths.
Their images under $\phi$ can be written:
\begin{eqnarray}
\label{tr1}
\phi([a,b])&=&[\phi(a),\phi(b)]=[b,c]\\
\label{tr2}
\phi([b,c])&=&[c,c^pa^{-r}b^{-q}]=c^pa^{-r}[a^r,c](b^{-q}[b^q,c]b^q)a^rc^{-p}\\
\phi([a,c])&=&[b,c^pa^{-r}b^{-q}]=[b,c^p]c^pa^{-r}[a^r,b]a^rc^{-p}
\label{tr3}
\end{eqnarray}

Using the identities 
$$[x^n,y]=(x^{n-1}[x,y]x^{1-n})\cdots(x[x,y]x^{-1})[x,y]$$ and
$$[x,y^n]=[x,y](y[x,y]y^{-1})\ldots(y^{n-1}[x,y]y^{1-n})$$ in the above,
we see that each of $\phi([a,b]),\phi([b,c])$ and $\phi([a,c])$ 
can be written as a product of conjugates
of $[a,b],[b,c],$ and $[a,c]$.
Interpreting this in terms of paths in $\R^2$, 
the closed paths $f\phi([a,b]),f\phi([b,c]),f\phi([a,c])$ can each be 
tiled by translates of copies of the parallelograms $f([a,b]),f([b,c]),f([a,c])$.

Geometrically, the images $f([a,b]),f([b,c]),f([c,a])$ are ``archtiles''
for a self-similar tiling of the plane.
Let $S=\{[a,b],[b,c],[a,c]\}$. The images $f(\phi(x))$, where $x\in S$
can be tiled by translates of tiles in $S$: the exact translates can be calculated
from the formulas in (\ref{tr1}),(\ref{tr2}),(\ref{tr3}).

Take $a,b,c$ to be the vectors $1,\lambda,$ and $\lambda^2\in \C$ respectively.
The equation (\ref{pqr}) implies that the argument of $\lambda$ is less than
$\pi/2$ degrees, so that $f[a,b],f[b,c],f[a,c]$ each have the same
orientation. (In particular this implies that the tiling will be non-overlapping.)
Then $f\phi(a)=f(b)=\lambda f(a)$, $f\phi(b)=f(c)=\lambda f(b)$, and $f\phi(c)$ is a path
whose endpoint is at $\lambda f(c)$.

Now the closed paths $f(\phi^n(x))$ for $x\in S$
converge after rescaling\footnote{a rigorous
argument for convergence can be made along the lines of Dekking \cite{Dek1,Dek2}} to the 
boundaries of certain tiles $T_1,T_2,T_3$, and $\lambda T_i$ can be tiled
by exactly by translates of $T_1,T_2$ and $T_3$. This is sufficient to
make a self-similar tiling of the plane.

\begin{figure}[htbp]
\centerline{\psfig{figure=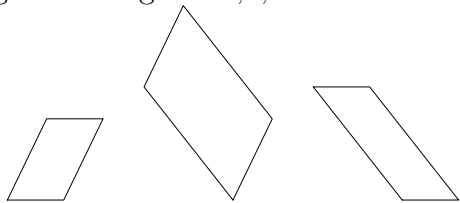,height=1in}}
\caption{\label{pars}The case $(p,q,r)=(1,2,1)$ ($\lambda\approx.696+1.436i$). The
three archtiles $f([a,b]),f([b,c]),f([a,c])$.}
\end{figure}

\begin{figure}[htbp]
\centerline{\psfig{figure=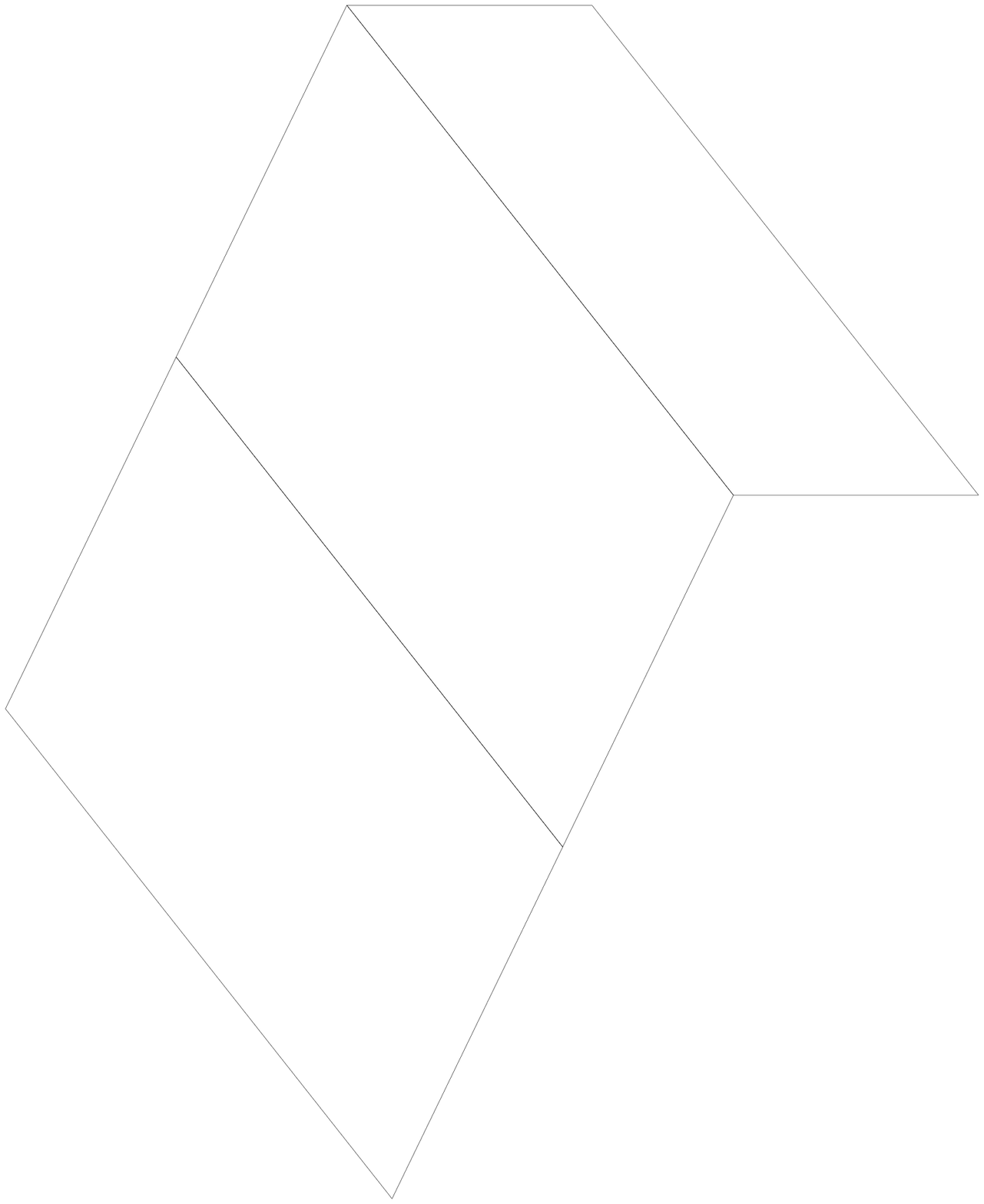,height=1in}
\psfig{figure=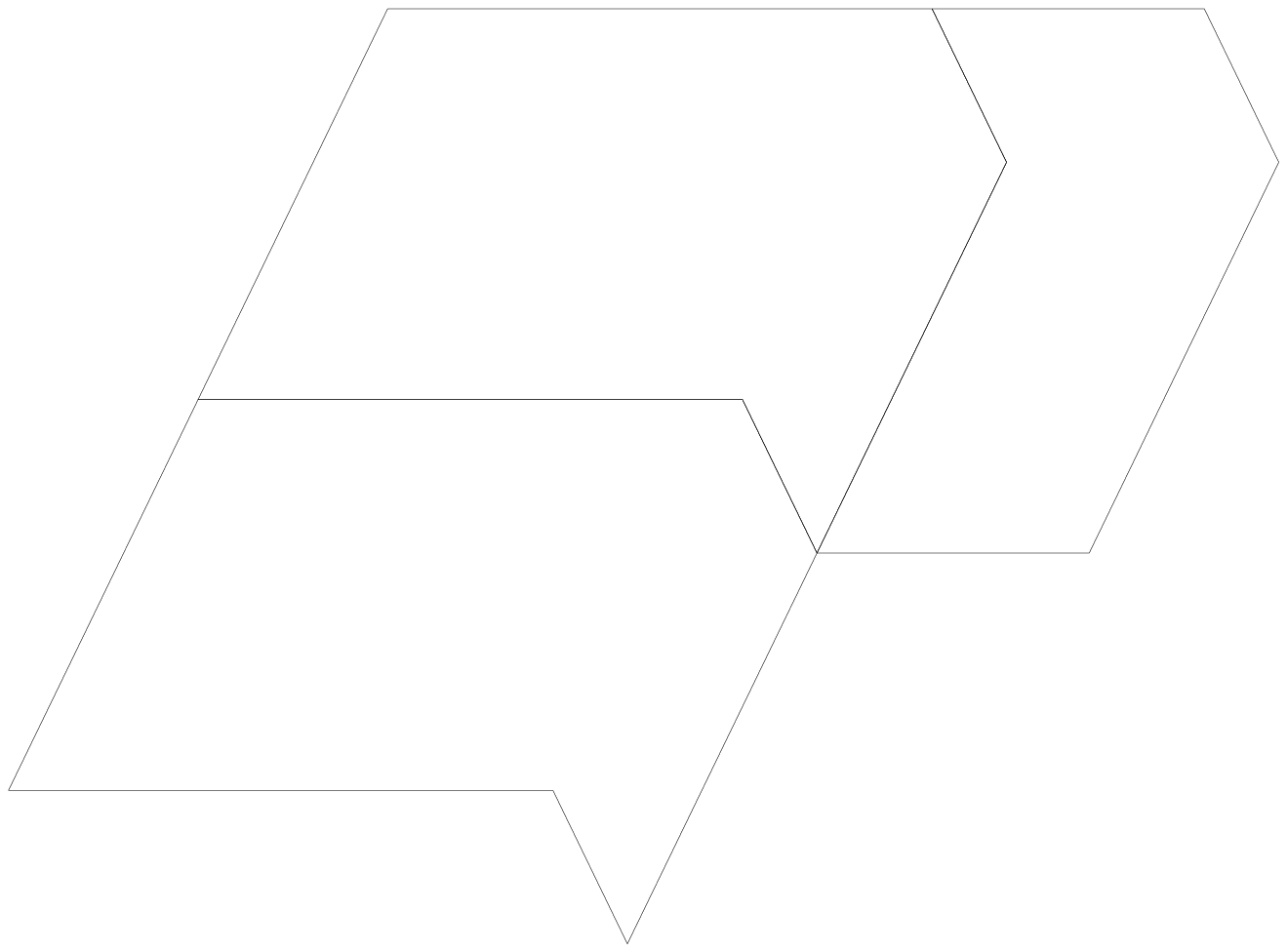,height=1in}
\psfig{figure=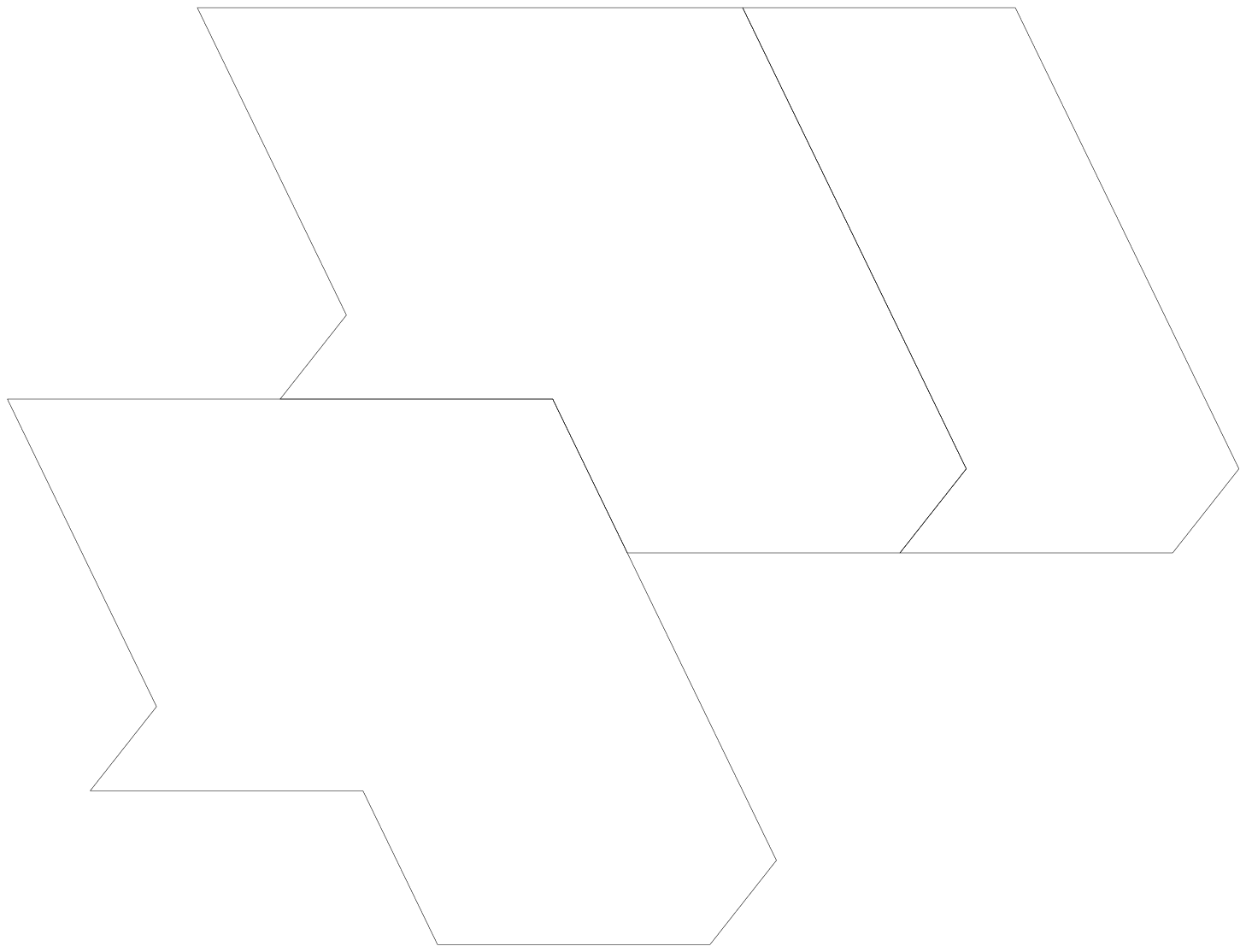,height=1in}
\psfig{figure=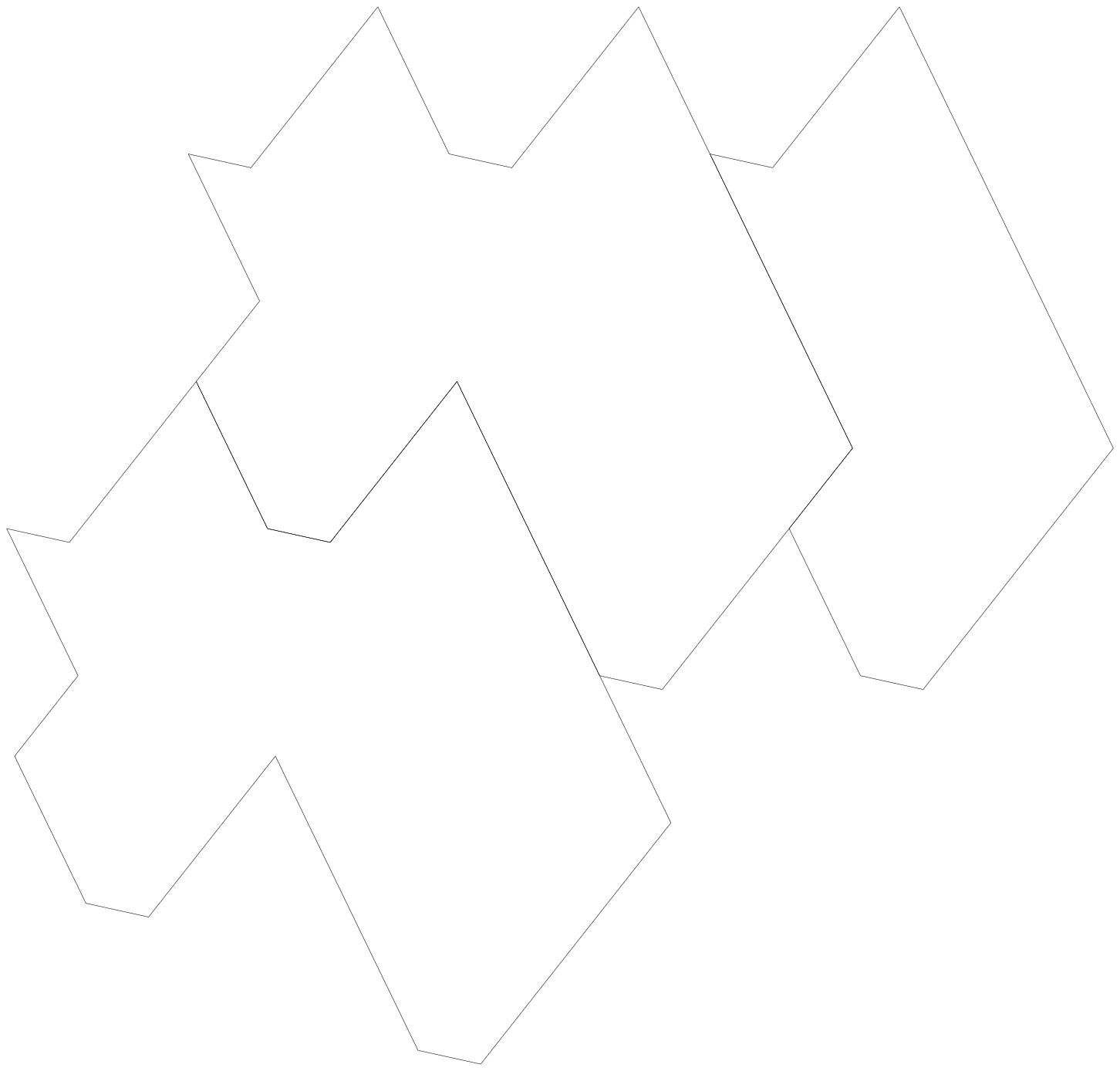,height=1in}
\psfig{figure=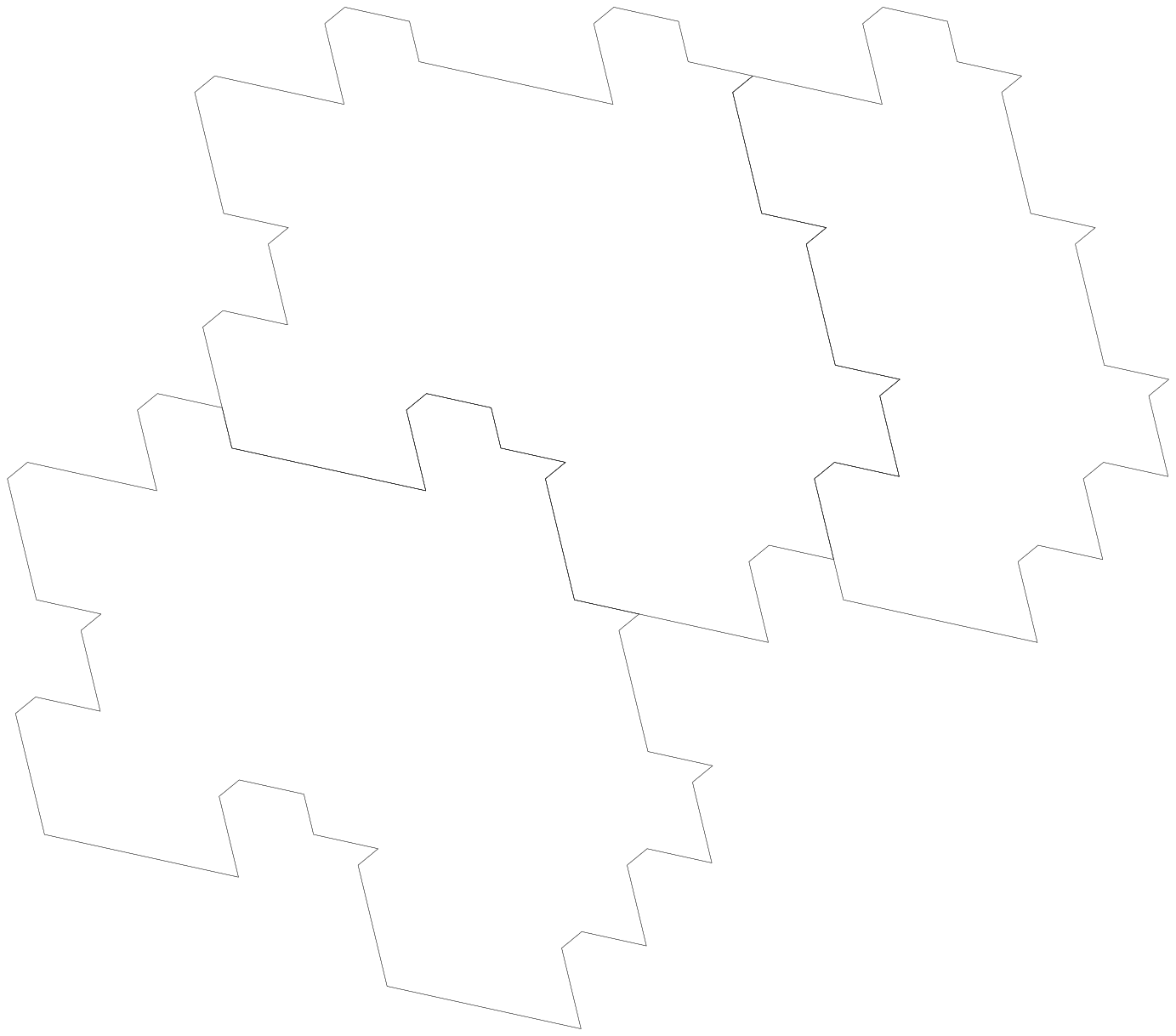,height=1in}
\psfig{figure=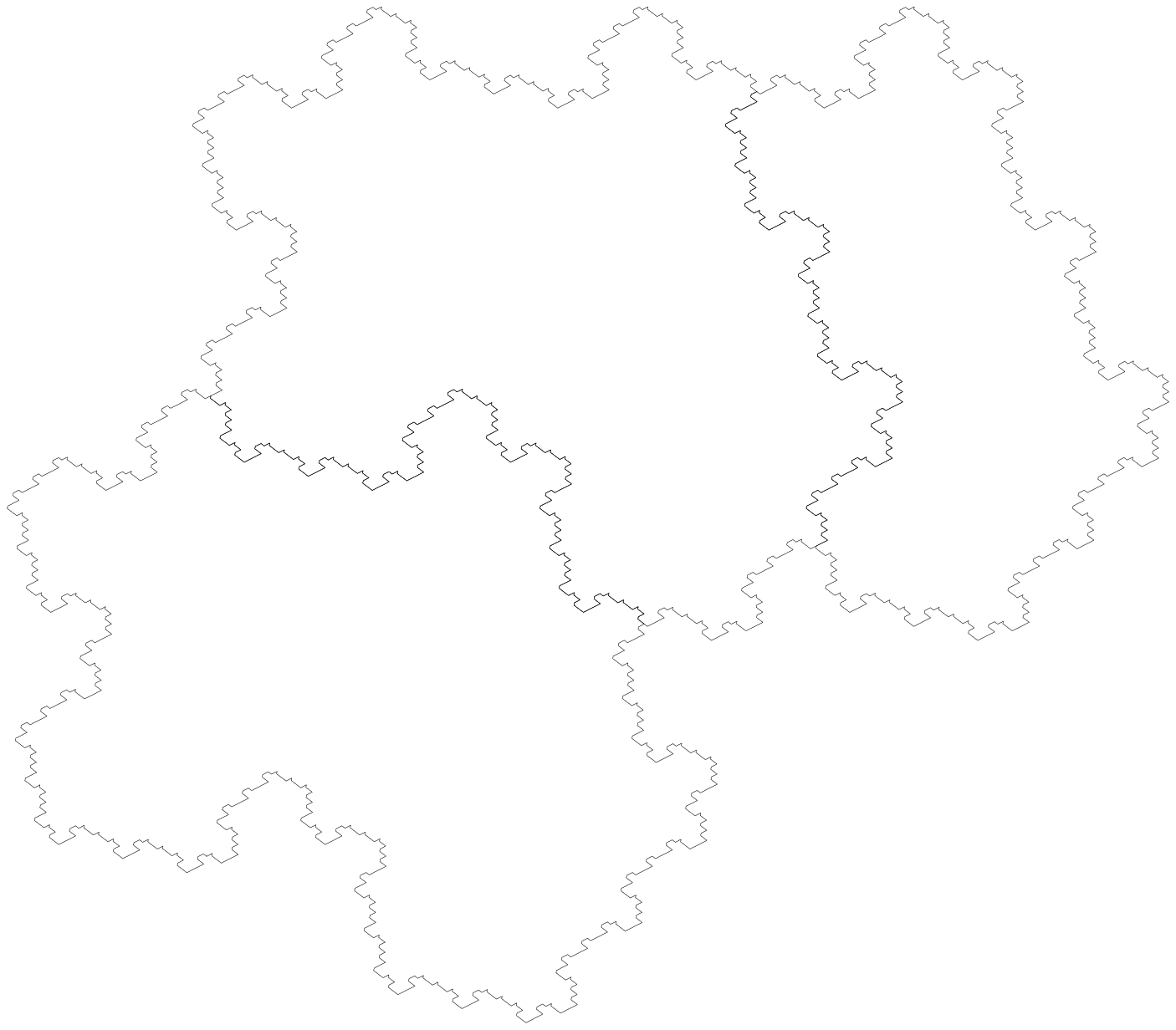,height=1in}}
\caption{\label{sst111}The successive approximations to the subdivision of $[b,c]$:
here we took $f\phi^n([b,c])$ for $n=2,3,4,5,6,$ and $10$.}
\end{figure}

\begin{figure}[htbp]
\centerline{\psfig{figure=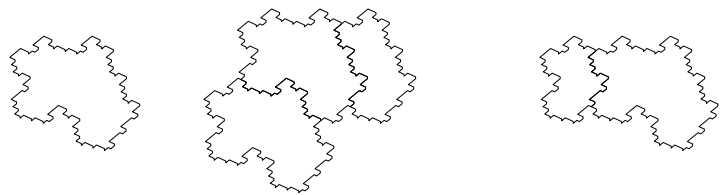,height=1.4in}}
\caption{\label{121tiles}The subdivisions of tiles $T_1,T_2,T_3$:
$T_1$ gives a $T_2$, $T_2$ gives two copies of $T_2$ and a copy of $T_3$, and
$T_3$ subdivides into a $T_2$ and a $T_1$.}
\end{figure}

An example is given in Figures \ref{pars},\ref{sst111},\ref{121tiles} and \ref{sst}.

For the case $n>3$, a similar construction will work.
Here we must use the free group on $n$ letters $a_1,\ldots,a_n$,
with an endomorphism $\phi(a_i)=a_{i+1}$ for $1\leq i\leq n-1$ and
$\phi(a_n)=a_n^{p}a_1^{-q}a_2^{-r}$.
The archtiles are the $(n^2-n)/2$ words/paths 
$[a_i,a_j]$ with $i<j$.

This same construction in fact
works with any endomorphism $\phi$ which has the property that
each $\phi([a_i,a_j])$ can be written as a product of conjugates of 
the $[a_{i'},a_{j'}]$. Since conjugates of the $[a_i,a_j]$ generate the
commutator subgroup $[F,F]$, this condition amounts to 
the requirement that the conjugates which appear have 
{\em nonnegative} exponent.

For an arbitrary endomorphism, one can also attempt to change the basis
of $F$ so that it has the correct form; a necessary condition is
that the induced linear map on $[F,F]/[F,[F,F]]\cong\Z^{(n^2-n)/2}$ 
have matrix with nonnegative coefficients.

\end{document}